\documentclass[14 pt, draft]{article}%[final]{dissert}
\usepackage[T2A]{fontenc}
\usepackage[cp1251]{inputenc}
\usepackage[ukrainian, english]{babel}
\usepackage[mathcal]{euscript}
\usepackage{amscd,amsfonts}
\usepackage{titlesec}
\usepackage{amsmath,amssymb,amsthm}
\usepackage{mathrsfs}
\usepackage{indentfirst}
\usepackage{graphicx}
\usepackage{amscd}

\graphicspath{}

\begin{document}
{\large

\begin{center}
 Pokutnyi O.A. 
\end{center}

\begin{center}
Bifurcation of solutions of the second order boundary value problems  in the Hilbert spaces.
\end{center}
{\small Conditions of the existence of solutions of linear and perturbed linear boundary value problems in the
Hilbert spaces for the second order evolution equation are obtained.
}

Consider the following boundary value problem (BVP) in the Hilbert spaces
\begin{equation} \label{eq:1}
y''(t, \varepsilon) + A(t)y(t, \varepsilon) = \varepsilon A_{1}(t)y(t, \varepsilon)  + f(t),
\end{equation}

\begin{equation} \label{eq:2}
 l(y(\cdot, \varepsilon), y'(\cdot, \varepsilon))^{T} = \alpha,
 \end{equation}
where
$y : J \rightarrow \mathcal{H}$ is a vector-function $y \in C^{2}(J, \mathcal{H})$, $J \subset \mathbb{R}$, the closed operator-valued function
$A(t)$ acts from $J$ into the dense domain $D = D(A(t)) \subset \mathcal{H}$
which is independent from $t$, $l$ is a linear and bounded operator which translates solutions of (\ref{eq:1}) into
the Hilbert space $\mathcal{H}_{1}$, $A_{1}(t)$ is a linear and bounded operator valued function $|||A_{1}||| = \sup_{t \in J}||A_{1}(t)|| < \infty$, $\alpha \in \mathcal{H}_{1}$.
%The main part of investigation is obtaining the necessary and sufficient conditions of the existence of solutions of nonlinear boundary value problem
%(\ref{eq:1}), (\ref{eq:2}) and linear boundary value problem (\ref{eq:3}).

{\bf Linear case.}  At first we find the necessary and sufficient conditions of the existence of solutions of linear nonhomogeneous boundary value problem \begin{equation} \label{eq:3}
 y''_{0}(t) + A(t)y_{0}(t) = f(t),~~~ l(y_{0}(\cdot), y'_{0}(\cdot))^{T} = \alpha.
\end{equation}Let $x_{1}(t) = y_{0}(t)$, $x_{2}(t) = y_{0}'(t)$, $x(t) = (x_{1}(t), x_{2}(t))^{T}$, then we can rewrite boundary value problem (\ref{eq:3}) in the form of the operator system
\begin{equation} \label{eq:4}
 x_{0}'(t) = B(t)x_{0}(t) + g(t),  lx_{0}(\cdot) = \alpha,
\end{equation}
where
\begin{equation}
 B(t) = \left( \begin{array}{ccc} 0 & I \\
                -A(t) & 0
               \end{array} \right),~~ g(t) = (0, f(t))^{T}.
\end{equation}

Denote by $U(t)$ the evolution operator of homogeneous equation $U'(t) = B(t)U(t), U(0) = I$. Then the set of solutions of (\ref{eq:4}) has the form
$$
x_{0}(t, c) = U(t)c + \int_{0}^{t}U(t)U^{-1}(\tau)g(\tau)d\tau.
$$
Substituting in $lx_{0}(\cdot) = \alpha$ we obtain the following operator equation
\begin{equation} \label{eq:5}
Qc = \alpha - l\int_{0}^{\cdot}U(\cdot)U^{-1}(\tau)g(\tau)d\tau, Q = lU(\cdot) : \mathcal{H} \rightarrow \mathcal{H}_{1}.
\end{equation}
Using the theory of strong generalized solutions \cite{2,3} we obtain the following result.

%\begin{theorem}
{\bf Theorem 1. } {\it 1. a) Boundary value problem (\ref{eq:3}) has strongly generalized solutions if and only if the following condition holds
\begin{equation} \label{eq:6}
\mathcal{P}_{N(\overline{Q}^{*})} \{ \alpha - l\int_{0}^{\cdot}U(\cdot)U^{-1}(\tau)g(\tau)d\tau \} = 0;
\end{equation}
if $\alpha - l\int_{0}^{\cdot}U(\cdot)U^{-1}(\tau)f(\tau)d\tau \in R(Q)$ then generalized solutions will be classical;

b) under condition (\ref{eq:6}) the set of solutions has the form
$$
x_{0}(t,c) = U(t)\mathcal{P}_{N(\overline{Q})}c + \overline{(G[g, \alpha])(t)},
$$
where $\mathcal{P}_{N(\overline{Q})}, \mathcal{P}_{N(\overline{Q}^{*})}$ are the orthoprojectors onto the kernel and cokernel of the operator
$\overline{Q}$ respectively,
$$
\overline{(G[g, \alpha])(t)} = \int_{0}^{t}U(t)U^{-1}(\tau)g(\tau)d\tau + \overline{Q}^{+}\{ \alpha - l\int_{0}^{\cdot}U(\cdot)U^{-1}(\tau)g(\tau)d\tau \}
$$
is a generalized Green's operator;

2. a) Boundary value problem (\ref{eq:3}) has strongly quasisolutions if and only if the following condition holds
\begin{equation} \label{eq:7}
\mathcal{P}_{N(\overline{Q}^{*})} \{ \alpha - l\int_{0}^{\cdot}U(\cdot)U^{-1}(\tau)f(\tau)d\tau \} \neq 0;
\end{equation}

b) under condition (\ref{eq:7}) the set of strongly quasisolutions has the form
$$
x_{0}(t,c) = U(t)\mathcal{P}_{N(\overline{Q})}c + \overline{(G[g, \alpha])(t)}.
$$
}

{\bf Bifurcation conditions.} a) Suppose that condition (\ref{eq:7}) is hold. We obtain the condition on $A_{1}(t)$ such that the perturbed boundary value problem
\begin{equation} \label{eq:1000}
x'(t,\varepsilon) = B(t)x(t, \varepsilon) + g(t) + \varepsilon B_{1}(t)x(t,\varepsilon),
\end{equation}
\begin{equation} \label{eq:1001}
lx(\cdot, \varepsilon) = \alpha,
\end{equation}
have the generalized solutions. Here is the operator-valued function $B_{1}(t)$ has the following form:
\begin{equation}
 B_{1}(t) = \left( \begin{array}{ccc} 0 & 0 \\
                 0 & A_{1}(t)
               \end{array} \right),~~ g(t) = (0, f(t))^{T},
\end{equation}
$x(t,\varepsilon) = (x_{1}(t,\varepsilon), x_{2}(t, \varepsilon))^{T}$, $x_{1}(t,\varepsilon) = y(t,\varepsilon), x_{2}(t, \varepsilon) = y'(t,\varepsilon)$. We will use the modification of the well-known Vishik-Lyusternik method. A solution of problem (\ref{eq:1000}), (\ref{eq:1001})
is sought in the form of a segment of the Laurent series in powers of the small parameter $\varepsilon$:
\begin{equation} \label{eq:1002}
x(t, \varepsilon) = \sum_{i = -1}^{+\infty}\varepsilon^{i}x_{i}(t) = \frac{x_{-1}(t)}{\varepsilon} + x_{0}(t) + \varepsilon x_{1}(t) + \varepsilon^{2}x_{2}(t) + ... .
\end{equation}
Substituting series (\ref{eq:1002}) into problem (\ref{eq:1000}), (\ref{eq:1001}) and equating the coefficients of $\varepsilon^{-1}$, we obtain the following boundary value problem for finding the coefficient $x_{-1}(t)$ of series (\ref{eq:1002}):
\begin{equation} \label{eq:1003}
x_{-1}^{'}(t) = B(t)x_{-1}(t),
\end{equation}
\begin{equation} \label{eq:1004}
lx_{-1}(\cdot) = 0.
\end{equation}
Problem (\ref{eq:1003}), (\ref{eq:1004}) has a family of solutions:
$$
x_{-1}(t, c_{-1}) = U(t)\mathcal{P}_{N(\overline{Q})}c_{-1}, c_{-1} \in \mathcal{H}.
$$
An arbitrary element $c_{-1}$  is determined by the condition for the solvability of the following  linear inhomogeneous boundary value problem for finding the coefficient $x_0(t)$ in series (\ref{eq:1002}):
\begin{equation} \label{eq:1005}
x_{0}^{'}(t) = B(t)x_{0}(t) + B_{1}(t)x_{-1}(t) + g(t),
\end{equation}
\begin{equation} \label{eq:1006}
lx_{0}(\cdot) = \alpha.
\end{equation}
A necessary and sufficient condition for the solvability of problem (\ref{eq:1005}), (\ref{eq:1006}) is given by
$$
\mathcal{P}_{N(\overline{Q}^{*})}\{ \alpha - l\int_{0}^{\cdot}U(\cdot)U^{-1}(\tau)(B_{1}(\tau)x_{-1}(\tau, c_{-1}) + g(\tau))d\tau \} = 0.
$$
 From this, in view of the form of $x_{-1}(t, c_{-1})$, we obtain an operator equation for
 $c_{-1} \in \mathcal{H}$:
\begin{equation} \label{eq:1007}
B_{0}c_{-1} =  \mathcal{P}_{N(\overline{Q}^{*})} \{ \alpha - l\int_{0}^{\cdot}U(\cdot)U^{-1}(\tau)g(\tau)d\tau\},
\end{equation}
where
$$
B_{0} = \mathcal{P}_{N(\overline{Q}^{*})}l\int_{0}^{\cdot}U(\cdot)U^{-1}(\tau)B_{1}(\tau)U(\tau)d\tau\mathcal{P}_{N(\overline{Q})}.
$$
A necessary and sufficient condition for the generalized solvability of this operator equation is
\begin{equation} \label{eq:1010}
\mathcal{P}_{N(\overline{B}_{0}^{*})}\mathcal{P}_{N(\overline{Q}^{*})} \{ \alpha - l\int_{0}^{\cdot}U(\cdot)U^{-1}(\tau)g(\tau)d\tau\} = 0.
\end{equation}
Suppose that $\mathcal{P}_{N(\overline{B}_{0}^{*})}\mathcal{P}_{N(\overline{Q}^{*})} = 0$. Then condition (\ref{eq:1010}) is hold. The solution set of operator equation for $c_{-1} \in \mathcal{H}$ has the form
$$
c_{-1} = \overline{c}_{-1} + \mathcal{P}_{N(\overline{B}_{0})}c_{\rho},  \forall c_{\rho} \in \mathcal{H},
$$
where
$$
\overline{c}_{-1} = \overline{B}_{0}^{+}\mathcal{P}_{N(\overline{Q}^{*})} \{ \alpha - l\int_{0}^{\cdot}U(\cdot)U^{-1}(\tau)g(\tau)d\tau\}.
$$
In view of the expression for $c_{-1}$, the homogeneous boundary value problem (\ref{eq:1005}), (\ref{eq:1006}) has a $\rho$ - parameter family of solutions
\begin{equation} \label{eq:1100}
x_{-1}(t, c_{\rho}) = \overline{x}_{-1}(t, \overline{c}_{-1}) + U(t)\mathcal{P}_{N(\overline{Q})}\mathcal{P}_{N(\overline{B}_{0})}c_{\rho},
\end{equation}
where
$$
\overline{x}_{-1}(t, \overline{c}_{-1}) = U(t)\mathcal{P}_{N(\overline{Q})}\overline{c}_{-1}.
$$
The general solution of problem (\ref{eq:1005}), (\ref{eq:1006}) has the form
$$
x_{0}(t, c_{0}) = U(t)\mathcal{P}_{N(\overline{Q})}c_{0} + F_{-1}(t) + K_{-1}(t)\mathcal{P}_{N(\overline{B}_{0})}c_{\rho},
$$
where
$$
F_{-1}(t) = \overline{(G[g + B_{1}\overline{x}_{-1}, \alpha])(t)}, K_{-1}(t) = \overline{(G[U, 0])(t)}\mathcal{P}_{N(\overline{Q})},
$$
where $c_0$ is an element of the space $\mathcal{H}$, which is determined at the next step from the condition for the solvability of the boundary value problem for finding the coefficient $x_{1}(t)$ in series (\ref{eq:1002}). To determine the coefficient $x_{1}(t)$ of $\varepsilon^{1}$ in series (\ref{eq:1002}), we obtain the following boundary value problem
\begin{equation}\label{eq:1111}
x_{1}^{'}(t) = B(t)x_{1}(t) + B_{1}(t)x_{0}(t, c_{0}),
\end{equation}
\begin{equation} \label{eq:1112}
lx_{1}(\cdot) = 0.
\end{equation}
Under condition of solvability
$$
\mathcal{P}_{N(\overline{Q}^{*})}l\int_{0}^{\cdot}U(\cdot)U^{-1}(\tau)B_{1}(\tau)x_{0}(\tau, c_{0})d\tau = 0,
$$
BVP (\ref{eq:1111}), (\ref{eq:1112}) has the set of solutions in the form
$$
x_{1}(t, c_{1}) = U(t)\mathcal{P}_{N(\overline{Q})}c_{1} + \overline{(G[B_{1}U\mathcal{P}_{N(\overline{Q})}c_{0} + F_{-1} + K_{-1}, 0])(t)}.
$$
The condition for the solvability of the boundary condition for the element $c_0$ is
\begin{equation} \label{eq:101010}
B_{0}c_{0} = -\mathcal{P}_{N(\overline{Q}^{*})}l(\int_{0}^{\cdot}U(\cdot)U^{-1}(\tau)B_{1}(\tau)F_{-1}(\tau)d\tau -  \end{equation}
$$
-\mathcal{P}_{N(\overline{Q}^{*})}l\int_{0}^{\cdot}U(\cdot)U^{-1}(\tau)B_{1}(\tau)K_{-1}(\tau)d\tau\mathcal{P}_{N(\overline{Q})}\mathcal{P}_{N(\overline{B}_{0})}c_{\rho}.
$$
From the condition $\mathcal{P}_{N(B_{0}^{*})}\mathcal{P}_{N(\overline{Q}^{*})} = 0$ follows solvability of equation (\ref{eq:101010}) with the set of solutions in the following form
$$
c_{0} = -B_{0}^{+}\mathcal{P}_{N(\overline{Q}^{*})}l\int_{0}^{\cdot}U(\cdot)U^{-1}(\tau)B_{1}(\tau)F_{-1}(\tau)d\tau -
$$
$$
-B_{0}^{+}\mathcal{P}_{N(\overline{Q}^{*})}l\int_{0}^{\cdot}U(\cdot)U^{-1}(\tau)B_{1}(\tau)K_{-1}(\tau)d\tau\mathcal{P}_{N(\overline{Q})}\mathcal{P}_{N(\overline{B}_{0})}c_{\rho} + \mathcal{P}_{N(\overline{B}_{0})}c_{\rho},
$$
$$
c_{0} = \overline{c}_{0} + D_{0}\mathcal{P}_{N(\overline{B}_{0})}c_{\rho}, ~~~\forall c_{\rho} \in \mathcal{H},
$$
where
$$
\overline{c}_{0} = -B_{0}^{+}\mathcal{P}_{N(\overline{Q}^{*})}l\int_{0}^{\cdot}U(\cdot)U^{-1}(\tau)B_{1}(\tau)F_{-1}(\tau)d\tau,
$$
$$
D_{0} = I - B_{0}^{+}\mathcal{P}_{N(\overline{Q}^{*})}l\int_{0}^{\cdot}U(\cdot)U^{-1}(\tau)B_{1}(\tau)K_{-1}(\tau)d\tau\mathcal{P}_{N(\overline{Q})}.
$$
Thus, problem (\ref{eq:1005}), (\ref{eq:1006})  has a $\rho$-parameter family of solutions:
$$
x_{0}(t, c_{0}) = \overline{x}_{0}(t, \overline{c}_{0}) + \overline{X}_{0}(t)\mathcal{P}_{N(\overline{B}_{0})}c_{\rho},~~~\forall c_{\rho} \in \mathcal{H},
$$
where
$$
\overline{x}_{0}(t, \overline{c}_{0}) = U(t)\mathcal{P}_{N(\overline{Q})}\overline{c}_{0} + F_{-1}(t),
$$
$$
\overline{X}_{0}(t) = U(t)\mathcal{P}_{N(\overline{Q})}D_{0} + K_{-1}(t).
$$
Then problem (\ref{eq:1111}), (\ref{eq:1112}) has a $\rho$-parameter family of solutions
$$
x_{0}(t, c_{\rho}) = \overline{x}_{0}(t, \overline{c}_{0}) + \overline{X}_{0}(t)\mathcal{P}_{N(\overline{B}_{0})}c_{\rho},
$$
where
$$
\overline{x}_{0}(t, \overline{c}_{0}) = U(t)\mathcal{P}_{N(\overline{Q})}\overline{c}_{0} + F_{-1}(t),
$$
$$
\overline{X}_{0}(t) = U(t)\mathcal{P}_{N(\overline{Q})}D_{0} + K_{-1}(t).
$$
Then problem (\ref{eq:1111}), (\ref{eq:1112}) has a $\rho$ - parameter family of solutions
$$
x_{1}(t, c_{1}) = U(t)\mathcal{P}_{N(\overline{Q})}c_{1} + F_{0}(t) + K_{0}(t)\mathcal{P}_{N(B_{0})}c_{\rho},
$$
where
$$
F_{0}(t) = \overline{(G[B_{1}U\mathcal{P}_{N(\overline{Q})}\overline{c}_{0} + F_{-1} + K_{-1}, 0])(t)},
$$
$$
K_{0}(t) = \overline{(G[B_{1}U\mathcal{P}_{N(\overline{Q})}D_{0}, 0])(t)},
$$
where $c_1$ is an element of the Hilbert space $\mathcal{H}$, which is determined at the next step from the condition for the solvability of the boundary value problem for finding the coefficient $x_2(t)$ in series (\ref{eq:1002}). By induction, we can prove that the coefficients $x_i(t)$ in series (\ref{eq:1002}) are determined by solving the boundary value problem
\begin{equation} \label{eq:1120}
x_{i}^{'}(t) = B(t)x_{i}(t) + B_{1}(t)x_{i - 1}(t, c_{i -1}),
\end{equation}
\begin{equation}\label{eq:1121}
lx_{i}(\cdot) = 0
\end{equation}
which  under condition of solvability has a $\rho$-parameter family of solutions
\begin{equation} \label{eq:1122}
x_{i}(t, c_{i}) = \overline{x}_{i}(t, \overline{c}_{i}) + \overline{X}_{i}(t)\mathcal{P}_{N(\overline{B}_{0})}c_{\rho},~~~\forall c_{\rho} \in \mathcal{H}
\end{equation}
where all the terms are determined by the iterative procedure
\begin{equation} \label{eq:1123}
\overline{x}_{i}(t, \overline{c}_{i}) = U(t)\mathcal{P}_{N(\overline{Q})}\overline{c}_{i} + F_{i -1}(t),
\end{equation}
\begin{equation} \label{eq:1124}
\overline{X}_{i}(t) = U(t)\mathcal{P}_{N(\overline{Q})}D_{i} + K_{i - 1}(t),
\end{equation}
\begin{equation} \label{eq:1125}
D_{i} = I - \overline{B}_{0}^{+}\mathcal{P}_{N(\overline{Q}^{*})}l\int_{0}^{\cdot}U(\cdot)U^{-1}(\tau)B_{1}(\tau)K_{i - 1}(\tau)d\tau\mathcal{P}_{N(\overline{Q})},
\end{equation}
\begin{equation} \label{eq:1126}
F_{i -1}(t) = \overline{(G[B_{1}U\mathcal{P}_{N(\overline{Q})}\overline{c}_{i - 1} + F_{i - 2} + K_{i - 2}, 0])(t)},
\end{equation}
\begin{equation} \label{eq:1127}
K_{i - 1}(t) = \overline{(G[B_{1}U\mathcal{P}_{N(\overline{Q})}D_{i - 1}, 0])(t)}.
\end{equation}
The convergence of series (\ref{eq:1002}) is proved in the same manner as in [11]. Thus, the following result holds.

{\bf Theorem 1.} {\it  The boundary value problem (\ref{eq:1000}), (\ref{eq:1001}) with the condition $\mathcal{P}_{N(\overline{B}_{0}^{*})}\mathcal{P}_{N(\overline{Q}^{*})} = 0$
 has a $\rho$-parameter family of  solutions in the form of the Laurent series segment
$$
x(t, c_{\rho}) = \sum_{i = -1}^{+ \infty} \varepsilon^{i} [\overline{x}_{i}(t, \overline{c}_{i}) + \overline{X}_{i}(t)\mathcal{P}_{N(\overline{B}_{0})}c_{\rho}], ~~~~\forall c_{\rho} \in \mathcal{H},
$$
whose coefficients are given by formulas (\ref{eq:1123})-(\ref{eq:1127}).
}

b) Suppose that condition (\ref{eq:6}) is hold. We obtain the condition on $A_{1}(t)$ such that the perturbed boundary value problem
\begin{equation} \label{eq:1000100}
x'(t,\varepsilon) = B(t)x(t, \varepsilon) + g(t) + \varepsilon B_{1}(t)x(t,\varepsilon),
\end{equation}
\begin{equation} \label{eq:1001100}
lx(\cdot, \varepsilon) = \alpha,
\end{equation}
have the generalized solutions.  A solution of problem (\ref{eq:1000100}), (\ref{eq:1001100})
is sought in the form of a segment of the Taylor series in powers of the small parameter $\varepsilon$:
\begin{equation} \label{eq:1002100}
x(t, \varepsilon) = \sum_{i = 0}^{+\infty}\varepsilon^{i}x_{i}(t) = x_{0}(t) + \varepsilon x_{1}(t) + \varepsilon^{2}x_{2}(t) + ... .
\end{equation}
Substituting series (\ref{eq:1002100}) into problem (\ref{eq:1000100}), (\ref{eq:1001100}) and equating the coefficients of $\varepsilon^{0}$, we obtain the following boundary value problem for finding the coefficient $x_{0}(t)$ of series (\ref{eq:1002100}):
\begin{equation} \label{eq:1003100}
x_{0}^{'}(t) = B(t)x_{0}(t) + g(t),
\end{equation}
\begin{equation} \label{eq:1004100}
lx_{0}(\cdot) = \alpha.
\end{equation}
Problem (\ref{eq:1003100}), (\ref{eq:1004100}) has a family of solutions:
$$
x_{0}(t,c) = U(t)\mathcal{P}_{N(\overline{Q})}c_{0} + \overline{(G[g, \alpha])(t)}.
$$
An arbitrary element $c_{0}$  is determined by the condition for the solvability of the following  linear inhomogeneous boundary value problem for finding the coefficient $x_1(t)$ in series (\ref{eq:1002100}):
\begin{equation} \label{eq:1005100}
x_{1}^{'}(t) = B(t)x_{1}(t) + B_{1}(t)x_{0}(t, c_{0}),
\end{equation}
\begin{equation} \label{eq:1006100}
lx_{1}(\cdot) = 0.
\end{equation}
A necessary and sufficient condition for the solvability of problem (\ref{eq:1005100}), (\ref{eq:1006100}) is given by
$$
\mathcal{P}_{N(\overline{Q}^{*})}\{ l\int_{0}^{\cdot}U(\cdot)U^{-1}(\tau)B_{1}(\tau)x_{0}(\tau, c_{0})d\tau \} = 0.
$$
 From this, in view of the form of $x_{0}(t, c_{0})$, we obtain an operator equation for
 $c_{0} \in \mathcal{H}$:
\begin{equation} \label{eq:1007100}
B_{0}c_{0} =  -\mathcal{P}_{N(\overline{Q}^{*})} l\int_{0}^{\cdot}U(\cdot)U^{-1}(\tau)B_{1}(\tau)\overline{(G[g, \alpha])(\tau)}d\tau.
\end{equation}
Under condition $\mathcal{P}_{N(\overline{B}_{0}^{*})}\mathcal{P}_{N(\overline{Q}^{*})} = 0$ the equation (\ref{eq:1007100}) is solvable. The solution set of operator equation for $c_{-1} \in \mathcal{H}$ has the form
$$
c_{0} = \overline{c}_{0} + \mathcal{P}_{N(\overline{B}_{0})}c_{\rho},  \forall c_{\rho} \in \mathcal{H},
$$
where
$$
\overline{c}_{0} = - \overline{B}_{0}^{+}\mathcal{P}_{N(\overline{Q}^{*})} l\int_{0}^{\cdot}U(\cdot)U^{-1}(\tau)B_{1}(\tau)\overline{(G[g, \alpha])(\tau)}d\tau.
$$
In view of the expression for $c_{0}$, the homogeneous boundary value problem (\ref{eq:1005100}), (\ref{eq:1006100}) has a $\rho$ - parameter family of solutions
\begin{equation} \label{eq:1100100}
x_{0}(t, c_{\rho}) = \overline{x}_{0}(t, \overline{c}_{0}) + U(t)\mathcal{P}_{N(\overline{Q})}\mathcal{P}_{N(\overline{B}_{0})}c_{\rho},
\end{equation}
where
$$
\overline{x}_{0}(t, \overline{c}_{0}) = U(t)\mathcal{P}_{N(\overline{Q})}\overline{c}_{0} + \overline{(G[g, \alpha])(t)}.
$$
The general solution of problem (\ref{eq:1005100}), (\ref{eq:1006100}) has the form
$$
x_{0}(t, c_{0}) = U(t)\mathcal{P}_{N(\overline{Q})}c_{0} + F_{-1}(t) + K_{-1}(t)\mathcal{P}_{N(\overline{B}_{0})}c_{\rho},
$$
where
$$
F_{-1}(t) = \overline{(G[g + B_{1}\overline{x}_{-1}, \alpha])(t)}, K_{-1}(t) = \overline{(G[U, 0])(t)}\mathcal{P}_{N(\overline{Q})},
$$
where $c_0$ is an element of the space $\mathcal{H}$, which is determined at the next step from the condition for the solvability of the boundary value problem for finding the coefficient $x_{1}(t)$ in series (\ref{eq:1002100}). To determine the coefficient $x_{1}(t)$ of $\varepsilon^{1}$ in series (\ref{eq:1002100}), we obtain the following boundary value problem
\begin{equation}\label{eq:1111100}
x_{1}^{'}(t) = B(t)x_{1}(t) + B_{1}(t)x_{0}(t, c_{0}),
\end{equation}
\begin{equation} \label{eq:1112100}
lx_{1}(\cdot) = 0.
\end{equation}
Under condition of solvability
$$
\mathcal{P}_{N(\overline{Q}^{*})}l\int_{0}^{\cdot}U(\cdot)U^{-1}(\tau)B_{1}(\tau)x_{0}(\tau, c_{0})d\tau = 0,
$$
BVP (\ref{eq:1111100}), (\ref{eq:1112100}) has the set of solutions in the form
$$
x_{1}(t, c_{1}) = U(t)\mathcal{P}_{N(\overline{Q})}c_{1} + \overline{(G[B_{1}U\mathcal{P}_{N(\overline{Q})}c_{0} + F_{-1} + K_{-1}, 0])(t)}.
$$
The condition for the solvability of the boundary condition for the element $c_0$ is
\begin{equation} \label{eq:101010100}
B_{0}c_{0} = -\mathcal{P}_{N(\overline{Q}^{*})}l(\int_{0}^{\cdot}U(\cdot)U^{-1}(\tau)B_{1}(\tau)F_{-1}(\tau)d\tau -  \end{equation}
$$
-\mathcal{P}_{N(\overline{Q}^{*})}l\int_{0}^{\cdot}U(\cdot)U^{-1}(\tau)B_{1}(\tau)K_{-1}(\tau)d\tau\mathcal{P}_{N(\overline{Q})}\mathcal{P}_{N(\overline{B}_{0})}c_{\rho}.
$$
From the condition $\mathcal{P}_{N(B_{0}^{*})}\mathcal{P}_{N(\overline{Q}^{*})} = 0$ follows solvability of equation (\ref{eq:101010}) with the set of solutions in the following form
$$
c_{0} = -B_{0}^{+}\mathcal{P}_{N(\overline{Q}^{*})}l\int_{0}^{\cdot}U(\cdot)U^{-1}(\tau)B_{1}(\tau)F_{-1}(\tau)d\tau -
$$
$$
-B_{0}^{+}\mathcal{P}_{N(\overline{Q}^{*})}l\int_{0}^{\cdot}U(\cdot)U^{-1}(\tau)B_{1}(\tau)K_{-1}(\tau)d\tau\mathcal{P}_{N(\overline{Q})}\mathcal{P}_{N(\overline{B}_{0})}c_{\rho} + \mathcal{P}_{N(\overline{B}_{0})}c_{\rho},
$$
$$
c_{0} = \overline{c}_{0} + D_{0}\mathcal{P}_{N(\overline{B}_{0})}c_{\rho}, ~~~\forall c_{\rho} \in \mathcal{H},
$$
where
$$
\overline{c}_{0} = -B_{0}^{+}\mathcal{P}_{N(\overline{Q}^{*})}l\int_{0}^{\cdot}U(\cdot)U^{-1}(\tau)B_{1}(\tau)F_{-1}(\tau)d\tau,
$$
$$
D_{0} = I - B_{0}^{+}\mathcal{P}_{N(\overline{Q}^{*})}l\int_{0}^{\cdot}U(\cdot)U^{-1}(\tau)B_{1}(\tau)K_{-1}(\tau)d\tau\mathcal{P}_{N(\overline{Q})}.
$$
Thus, problem (\ref{eq:1005100}), (\ref{eq:1006100})  has a $\rho$-parameter family of solutions:
$$
x_{0}(t, c_{0}) = \overline{x}_{0}(t, \overline{c}_{0}) + \overline{X}_{0}(t)\mathcal{P}_{N(\overline{B}_{0})}c_{\rho},~~~\forall c_{\rho} \in \mathcal{H},
$$
where
$$
\overline{x}_{0}(t, \overline{c}_{0}) = U(t)\mathcal{P}_{N(\overline{Q})}\overline{c}_{0} + F_{-1}(t),
$$
$$
\overline{X}_{0}(t) = U(t)\mathcal{P}_{N(\overline{Q})}D_{0} + K_{-1}(t).
$$
Then problem (\ref{eq:1111100}), (\ref{eq:1112100}) has a $\rho$-parameter family of solutions
$$
x_{0}(t, c_{\rho}) = \overline{x}_{0}(t, \overline{c}_{0}) + \overline{X}_{0}(t)\mathcal{P}_{N(\overline{B}_{0})}c_{\rho},
$$
where
$$
\overline{x}_{0}(t, \overline{c}_{0}) = U(t)\mathcal{P}_{N(\overline{Q})}\overline{c}_{0} + F_{-1}(t),
$$
$$
\overline{X}_{0}(t) = U(t)\mathcal{P}_{N(\overline{Q})}D_{0} + K_{-1}(t).
$$
Then problem (\ref{eq:1111100}), (\ref{eq:1112100}) has a $\rho$ - parameter family of solutions
$$
x_{1}(t, c_{1}) = U(t)\mathcal{P}_{N(\overline{Q})}c_{1} + F_{0}(t) + K_{0}(t)\mathcal{P}_{N(B_{0})}c_{\rho},
$$
where
$$
F_{0}(t) = \overline{(G[B_{1}U\mathcal{P}_{N(\overline{Q})}\overline{c}_{0} + F_{-1} + K_{-1}, 0])(t)},
$$
$$
K_{0}(t) = \overline{(G[B_{1}U\mathcal{P}_{N(\overline{Q})}D_{0}, 0])(t)},
$$
where $c_1$ is an element of the Hilbert space $\mathcal{H}$, which is determined at the next step from the condition for the solvability of the boundary value problem for finding the coefficient $x_2(t)$ in series (\ref{eq:1002100}). By induction, we can prove that the coefficients $x_i(t)$ in series (\ref{eq:1002100}) are determined by solving the boundary value problem
\begin{equation} \label{eq:1120100}
x_{i}^{'}(t) = B(t)x_{i}(t) + B_{1}(t)x_{i - 1}(t, c_{i -1}),
\end{equation}
\begin{equation}\label{eq:1121100}
lx_{i}(\cdot) = 0
\end{equation}
which  under condition of solvability has a $\rho$-parameter family of solutions
\begin{equation} \label{eq:1122100}
x_{i}(t, c_{i}) = \overline{x}_{i}(t, \overline{c}_{i}) + \overline{X}_{i}(t)\mathcal{P}_{N(\overline{B}_{0})}c_{\rho},~~~\forall c_{\rho} \in \mathcal{H}
\end{equation}
where all the terms are determined by the iterative procedure
\begin{equation} \label{eq:1123100}
\overline{x}_{i}(t, \overline{c}_{i}) = U(t)\mathcal{P}_{N(\overline{Q})}\overline{c}_{i} + F_{i -1}(t),
\end{equation}
\begin{equation} \label{eq:1124100}
\overline{X}_{i}(t) = U(t)\mathcal{P}_{N(\overline{Q})}D_{i} + K_{i - 1}(t),
\end{equation}
\begin{equation} \label{eq:1125100}
D_{i} = I - \overline{B}_{0}^{+}\mathcal{P}_{N(\overline{Q}^{*})}l\int_{0}^{\cdot}U(\cdot)U^{-1}(\tau)B_{1}(\tau)K_{i - 1}(\tau)d\tau\mathcal{P}_{N(\overline{Q})},
\end{equation}
\begin{equation} \label{eq:1126100}
F_{i -1}(t) = \overline{(G[B_{1}U\mathcal{P}_{N(\overline{Q})}\overline{c}_{i - 1} + F_{i - 2} + K_{i - 2}, 0])(t)},
\end{equation}
\begin{equation} \label{eq:1127100}
K_{i - 1}(t) = \overline{(G[B_{1}U\mathcal{P}_{N(\overline{Q})}D_{i - 1}, 0])(t)}.
\end{equation}
The convergence of series (\ref{eq:1002100}) is proved in the same manner as in [11]. Thus, the following result holds.

{\bf Theorem 2.} {\it  The boundary value problem (\ref{eq:1000100}), (\ref{eq:1001100}) with the condition $\mathcal{P}_{N(\overline{B}_{0}^{*})}\mathcal{P}_{N(\overline{Q}^{*})} = 0$
 has a $\rho$-parameter family of  solutions in the form of the Laurent series segment
$$
x(t, c_{\rho}) = \sum_{i = -1}^{+ \infty} [\overline{x}_{i}(t, \overline{c}_{i}) + \overline{X}_{i}(t)\mathcal{P}_{N(\overline{B}_{0})}c_{\rho}], ~~~~\forall c_{\rho} \in \mathcal{H},
$$
whose coefficients are given by formulas (\ref{eq:1123100})-(\ref{eq:1127100}).
}

}
\end{document}